\documentclass[11pt,english]{article}

\usepackage{amsfonts,amsmath,amsthm,amscd,amssymb,latexsym,cite,verbatim}
\usepackage[ukrainian,russian,english]{babel}

\theoremstyle{plain}
\newtheorem{theorem}{Theorem}[section]
\newtheorem{lemma}{Lemma}[section]

\newtheorem{corollary}{Corollary}[section]

\newtheorem{definition}{Definition}[section]
\theoremstyle{definition}

\newtheorem{remark}{Remark}[section]
\newcommand{\keywords}{\textbf{Key words. }\medskip}
\newcommand{\subjclass}{\textbf{MSC 2010. }\medskip}
\renewcommand{\abstract}{\textbf{Abstract. }\medskip}

\numberwithin{equation}{section}
\begin{document}

\title{Локализованные режимы с обострением для квазилинейных дважды вырождающихся
параболических уравнений}

\author{А.\,Е.~Шишков, Е.\,А.~Евгеньева}

\date{}

\maketitle

\begin{abstract}
Изучаются режимы с сингулярным обострением для широкого класса квазилинейных
параболических уравнений второго порядка. На основе энергетических методов
устанавливаются в определенном смысле точные оценки финального профиля обобщенного
решения в окрестности времени обострения в зависимости от скорости нарастания глобальной
энергии этого решения.
\end{abstract}



\section{\bf Введение и формулировка основных результатов}

В ограниченной цилиндрической области $Q=(0,T)\times\Omega$, $1\leqslant T<\infty$,
$\Omega\subset\mathbb{R}^n$, $n\geqslant1$, $\partial\Omega\in C^2$, рассматривается
уравнение
\begin{equation}\label{equat}
(|u|^{q-1}u)_t-\sum_{i=1}^n(a_i(t,x,u,\nabla u))_{x_i}=0,\quad q=\text{const}>0.
\end{equation}
Здесь функции $a_i(t,x,s,\xi)$, $i=1,2,...,n$ являются непрерывными функциями всех своих
аргументов и удовлетворяют следующим условиям коэрцитивности и роста:
\begin{equation}\label{coercit}
d_0|\xi|^{p+1}\leqslant\sum_{i=1}^na_i(t,x,s,\xi)\xi_i\quad
\forall(t,x,s,\xi)\in\bar{Q}\times\mathbb{R}^1\times\mathbb{R}^n,\
d_0=\textrm{const}>0,\ p=const>q,
\end{equation}
\begin{equation}\label{growth}
|a_i(t,x,s,\xi)|\leqslant d_1|\xi|^p\quad
\forall(t,x,s,\xi)\in\bar{Q}\times\mathbb{R}^1\times\mathbb{R}^n,\ i=1, ..., n,\
d_1=\textrm{const}<\infty.
\end{equation}
Для произвольной функции $u_0\in L_{q+1}(\Omega)$ обозначим через $U_{u_0}$ множество
всех слабых (энергетических) решений уравнения \eqref{equat}, удовлетворяющих начальному
условию
\begin{equation}\label{initial}
u(0,x)=u_0(x)\quad\forall x\in\Omega.
\end{equation}

\begin{definition}\label{Def1}
Функция $u(t,x)\in C_{loc}([0,T);L_{q+1}(\Omega))$ называется слабым \\(энергетическим)
решением задачи \eqref{equat}, \eqref{initial}, если
\begin{enumerate}
\item[i)] $u(t,\cdot)\in L_{p+1,\ loc}([0,T);W^1_{p+1}(\Omega))$;

\item[ii)] $(|u(t,\cdot)|^{q-1}u(t,\cdot))_t\in L_{\frac{p+1}p, \
loc}([0,T);(\stackrel{\circ}{W}^1_{p+1}(\Omega))^*)$;

\item[iii)] выполнено интегральное тождество:
\begin{equation}\label{ident1}
\int_0^\tau\langle(|u|^{q-1}u)_t,\eta\rangle dt+
\int_0^\tau\int_{\Omega}\sum_{i=1}^na_i(t,x,u,\nabla u)\eta_{x_i}dxdt=0
\end{equation}
для произвольной функции $\eta(t,\cdot)\in
L_{p+1}((0,\tau);\stackrel{\circ}{W}^1_{p+1}(\Omega))$ с любым $\tau< T$;

\item[iv)] выполнено начальное условие \eqref{initial}.
\end{enumerate}
\end{definition}
Пусть теперь $F(t)>0$ --- произвольная монотонно неубывающая на интервале $[0,T)$
функция такая, что $F(t)\rightarrow\infty$ при $t\rightarrow T$. Введем подмножество
$U_{u_0,F}\subset U_{u_0}$ всех энергетических решений $u$ задачи
\eqref{equat}--\eqref{initial}, удовлетворяющих оценке:
\begin{equation}\label{setU}
\begin{aligned}
&\sup_{0<\tau<t}h(\tau)+E(t)=\sup_{0<\tau<t}h_u(\tau)+E_u(t):=\sup_{0<\tau<t}\int_{\Omega}|u(\tau,x)|^{q+1}dx+
\\&+\int_0^t\int_{\Omega}|\nabla_xu(\tau,x)|^{p+1}dxd\tau\leqslant F(t)\quad
\forall\,t\in(0,T).
\end{aligned}
\end{equation}
Отметим, что  решения $u$ задачи \eqref{equat}--\eqref{initial}, энергетические функции
которых обладают ростом, определяемым оценкой \eqref{setU}, могут быть получены
различными способами. Например, при помощи граничных данных (граничного режима),
обостряющихся в момент времени $t=T$:
\begin{equation}\label{Dir1}
    u(t,x)\Big|_{\partial\Omega}=f(t,x)\rightarrow\infty\text{ при }t\rightarrow T
\end{equation}
либо
\begin{equation}\label{Neum1}
   \frac{\partial u(t,x)}{\partial N}\biggr|_{\partial\Omega}:=\sum_{i=1}^na_i(t,x,u,\nabla u)\ \nu_i
   =g(t,x)\rightarrow\infty\text{ при }t\rightarrow T,
\end{equation}
где $\nu=\nu(x)=(\nu_1, ..., \nu_n)$ --- единичный вектор внешней нормали к
$\partial\Omega$ в точке $x$. Несложно проверить (см. \cite{KSSh}), что энергетические
функции произвольного решения $u$ задачи \eqref{equat}, \eqref{initial}, \eqref{Dir1}
удовлетворяют оценке \eqref{setU} с
\begin{equation}\label{Fdef}
\begin{aligned}
&F(t):=F_f(t)=\sup_{0<\tau<t}\int_{\Omega}|\overline{f}(\tau,x)|^{q+1}dx+
\int_0^t\int_{\Omega}|\nabla_x\overline{f}(\tau,x)|^{p+1}dxd\tau+
\\&+\left(\int_0^t\left(\int_{\Omega}|\overline{f}(\tau,x)|^{q+1}dx\right)^\frac1{q+1}d\tau\right)^{q+1},
\end{aligned}
\end{equation}
где $\overline{f}(t,x)$ --- это любое продолжение функции ${f}(t,x)$ с
$[0,T)\times\partial\Omega$ на всю цилиндрическую область $Q$. Понятно, что точность
соответствующей оценки \eqref{setU} связана с оптимальностью выбора продолжения
$\overline{f}$.

Существование энергетических обобщенных решений задачи \eqref{equat}, \eqref{initial},
\eqref{Dir1} либо \eqref{Neum1} изучались многими авторами (в случае \eqref{Dir1} см.,
например, \cite{AL1}, \cite{Ivan.exist2}, \cite{Ivan.exist1} и имеющиеся там ссылки).
Нас интересует вопрос о поведении этих решений в окрестности времени обострения режима,
то есть при $t\rightarrow T$. В случае различных модельных уравнений типа нестационарной
ньютоновской или неньютоновской медленной или нейтральной диффузии (фильтрации),
допускающих построение автомодельных (или аппроксимативно автомодельных) решений и
теорем сравнения, такой анализ на основе соответствующей барьерной техники приводился
многими исследователями (см., например, \cite{SGKM1}, \cite{GH1}, \cite{CE}, \cite{GG},
\cite{Ven1}). Были описаны локализованные граничные режимы, то есть режимы, для которых
$\Omega\setminus\Omega_c\neq\emptyset$, где $\Omega_c$ --- это множество сингулярности
соответствующего решения $u$:
\begin{equation}\label{singSet}
\Omega_c:=\Omega_c(u)=\left\{x\in\Omega:\overline{\lim_{t\rightarrow
T}}\,u(t,x)=\infty\right\}.
\end{equation}
Найдены некоторые геометрические характеристики  множеств $\Omega_c$ в зависимости от
аналитической структуры соответствующих граничных режимов.

В \cite{ShSh1} был предложен и развит некоторый вариант метода локальных энергетических
оценок для изучения локализации граничных режимов с обострением для общих квазилинейных
параболических уравнений. Этот метод не использует никаких вариантов барьерной техники
или теорем сравнения и применим к очень широким классам уравнений, включая
параболические уравнения высоких порядков (см. \cite{GSh1}, \cite{GSh2}) и некоторые
псевдопараболические уравнения. В частности, в \cite{ShSh1} (см. также Теорема 6.3.1 из
\cite{KSSh}) доказано следующее утверждение. Если $u(t,x)$  --- произвольное
энергетическое решение задачи \eqref{equat}, \eqref{initial} из множества $U_{u_0,F}$
(см. определение \eqref{setU}) с функцией
\begin{equation}\label{S-reg}
F=F(t)\leqslant
F_0(t):=\omega(T-t)^{-\frac{q+1}{p-q}}\quad\forall\,t<T,\,\omega=const>0,
\end{equation}
то справедлива следующая равномерная по $t<T$ энергетическая априорная оценка:
\begin{equation}\label{S-est}
\begin{aligned}
&\sup_{0<\tau<t}h(\tau,s)+E(t,s)=\sup_{0<\tau<t}h_{u}(\tau,s)+E_{u}(t,s):=\sup_{0<\tau<t}\int_{\Omega(s)}|u(\tau,x)|^{q+1}dx+
\\&+\int_0^t\int_{\Omega(s)}|\nabla_xu(\tau,x)|^{p+1}dxd\tau<C<\infty\quad
\forall\,t<T,\quad\forall\,s>\delta(\omega),\,\delta(\omega)\rightarrow0\text{ при
}\omega\rightarrow0,
\end{aligned}
\end{equation}
где $\Omega(s)$ семейство подобластей, определенных соотношением:
\begin{equation}\label{Omega(s)}
\Omega(s):=\{x\in\Omega:d(x):=dist(x,\partial\Omega)>s\}.
\end{equation}
\begin{remark}\label{rem1}
В силу внутренних оценок максимума модулей энергетических решений уравнения
\eqref{equat} (см., например, \cite{Ivan.max1}) из оценки \eqref{S-est} вытекает
равномерная по $t<T$ ограниченность $|u(t,x)|$ для любой точки
$x\in\Omega(\delta(\omega))$. Следовательно в силу определения \eqref{singSet} имеем
\begin{equation*}
\Omega_c(u)\subset\Omega\setminus\Omega(\delta(\omega)),
\end{equation*}
а значит режим, определяемый в \eqref{S-reg} является локализованным. Более того,
дополнительный анализ показывает, что в общей ситуации этот $\Omega_c(u)\neq\emptyset$ и
в терминологии \cite{SGKM1} режим \eqref{S-reg} является локализованным $S$--режимом.
\end{remark}
\begin{remark}\label{rem2}
Если режим $F$ удовлетворяет асимптотически более сильной, чем \eqref{S-reg} оценке:
\begin{equation}\label{LS-reg}
F(t)\leqslant F_1(t):=\omega(t)(T-t)^{-\frac{q+1}{p-q}}\quad\forall\,t<T,
\end{equation}
где $\omega(t)$ --- некоторая неотрицательная невозрастающая функция такая, что
\\$\lim_{t\rightarrow T}\omega(t)=0$, то в силу оценки \eqref{S-est}
$\Omega_c(u)\cap\Omega=\emptyset$ и, значит, в терминологии \cite{SGKM1} режим
\eqref{LS-reg} является локализованным LS--режимом.
\end{remark}
Для LS--режимов существует функция
\begin{equation}\label{Profile1}
W(x):=\overline{\lim_{t\rightarrow T}}u(t,x)<\infty\quad\forall\,x\in\Omega,
\end{equation}
которую мы называем финальным профилем решения и которая является тонкой
характе-\\ристикой этого решения. Актуальным является получение точных зависящих от
$\omega(t)$ оценок сверху и снизу для этого финального профиля. Этой задаче посвящена
данная работа. Основным результатом является следующая теорема.

\begin{theorem}\label{Th1}
Пусть $u(t,x)$ --- произвольное энергетическое решение задачи \eqref{equat},
\eqref{initial} из множества $U_{u_0,F}$, где $F$ --- это локализованный LS-режим,
определяемый оценкой \eqref{LS-reg} с функцией
\begin{equation}\label{regLS}
\omega(t)=\omega_0(T-t)^{\beta} \quad\forall\,t<T,
\quad\omega_0=const>0,\,0<\beta=const<\beta_0=\frac{q+1}{p-q}-\frac1p.
\end{equation}
Тогда существует константа $G>0$ и значение $\hat{s}>0$, зависящие только от известных
параметров задачи, такие, что для решения $u$ справедлива следующая равномерная по
$t\leqslant T$ энергетическая оценка:
\begin{equation}\label{main}
\begin{aligned}
&E_u(t,s)+ \sup_{0<\tau<t}h_u(\tau,s)\leqslant G\omega_0^{\frac{q+1}{\beta(p-q)}}
s^{-\nu} \quad\quad\forall\,t\leqslant T,\,\forall\,s\in(0,\hat{s}).
%\\&\nu=\frac{(n(p-q)+(q+1)(p+1))(q+1-\beta(p-q))}{\beta(p-q)^2}
\end{aligned}
\end{equation}
где $\nu=\frac{(n(p-q)+(q+1)(p+1))(q+1-\beta(p-q))}{\beta(p-q)^2}$.
\end{theorem}

\begin{remark}\label{rem4}
Ограничение сверху на величину $\beta:\beta<\beta_0$ является чисто техническим. При
$\beta_0\leqslant\beta<\frac{q+1}{p-q}$ наш анализ требует дополнительных шагов, которые
мы приведем в следующей публикации. Также планируется перестроить анализ для случая
нестепенных граничных режимов, например при $\omega(t)=\left(-\ln(T-t)\right)^{-l}$,
$l>0$.
\end{remark}

\begin{corollary}\label{CorProfile}
В условиях теоремы \ref{Th1} для финального профиля $W(x)$ решения $u$ задачи
\eqref{equat}, \eqref{initial} при дополнительном условии $0<p-1<q<1$ выполняется
оценка:
\begin{equation}\label{main1}
W(x)\leqslant C_1\omega_0^{\frac1{\beta(p-q)}}
d(x)^{-\mu}\quad\forall\,x\in\Omega\setminus\Omega(2\hat{s}),
\end{equation}
где $C_1<\infty$ зависит только от известных параметров задачи, $\hat{s}$ из
\eqref{main},
\\а $\mu=\frac{n(p-q)+(p+1)(q+1)-\beta(p-q)(p+1)}{\beta(p-q)^2}$.
\end{corollary}

\section{\bf Финальный профиль решения}

{\bf Доказательство теоремы \ref{Th1}.} Из определения \eqref{S-est} следует, что
энергетические функции $E(t,s)$, $h(t,s)$ рассматриваемого решения $u$ являются
невозрастающими функциями аргумента $s$ при любом $t\leqslant T$. Кроме того в силу
замечания \ref{rem2} имеем
\begin{equation}\label{add1}
E(t,s)+\sup_{0<\tau<t}h(\tau,s)<\infty\quad\forall\,s>0,\,\forall\,t\leqslant T.
\end{equation}
Зафиксируем теперь числа
\begin{equation}\label{add2}
\xi\in(0,1);\quad\alpha_1\in(p^{-1},\alpha),\quad\alpha:=\frac{q+1}{p-q}-\beta.
\end{equation}
В силу \eqref{add1} имеется альтернатива: либо
\begin{equation}\label{add3}
E(T,{s})+\sup_{0<\tau<T}h(\tau,s)\leqslant2\omega_0T^{-\alpha_1}\xi^{-\alpha}\quad\forall\,s>0,
\end{equation}
либо существует такое значение $\bar{s}\in(0,s_\Omega)$, где $s_\Omega$ из
\eqref{s_Omega}, что
\begin{equation}\label{4.1}
E(T,{s})+\sup_{0<\tau<T}h(\tau,s)>2\omega_0T^{-\alpha_1}\xi^{-\alpha}\quad\forall\,s\in(0,\bar{s}).
\end{equation}
Начнем наш анализ с основного случая \eqref{4.1}. Для произвольной точки
$\tilde{s}\in(0,{\bar{s}})$ определим конечную возрастающую последовательность
$\{t_j\}=\{t_j(\tilde{s})\}$, $j=1,2,...$, $t_0=0$, при помощи непрерывной функции
$\Gamma_{\tilde{s}}(\cdot):[0,t']\rightarrow[t_1,T]$, определяемой следующим равенством:
\begin{equation}\label{11.3}
\big(\Gamma_{\tilde{s}}(t)-t\big)^{-\alpha}=\frac{\xi^\alpha}{\omega_0T^{\alpha-\alpha_1}}
\Big( E(\Gamma_{\tilde{s}}(t),\tilde{s})- E(t,\tilde{s})+
\sup_{t<\tau<\Gamma_{\tilde{s}}(t)} h(\tau,\tilde{s})\Big).
\end{equation}
Значение $t_1=t_1(\tilde{s})=\Gamma_{\tilde{s}}(0)$ определяется следующим равенством:
\begin{equation}\label{11.4}
t_1^{-\alpha}=\frac{\xi^\alpha}{\omega_0T^{\alpha-\alpha_1}}\left( E(t_1,\tilde{s})+
\sup_{0<\tau<t_1} h(\tau,\tilde{s})\right),
\end{equation}
а $t'$ определяется из соотношения:
\begin{equation}\label{11.5}
(T-t')^{-\alpha}=\frac{\xi^\alpha}{\omega_0T^{\alpha-\alpha_1}} \Big( E(T,\tilde{s})-
E(t',\tilde{s})+ \sup_{t'<\tau<T} h(\tau,\tilde{s})\Big).
\end{equation}
В силу определения \eqref{11.4} и предположения \eqref{4.1} имеем
\begin{equation}\label{11.61}
\begin{aligned}
&\left( E(t_1,\tilde{s})+\sup_{0<\tau<t_1} h(\tau,\tilde{s})\right)t_1^{\alpha}=
\frac{\omega_0T^{\alpha-\alpha_1}}{\xi^\alpha}<
\\&<\frac12 \left(
E(T,\tilde{s})+\sup_{0<\tau<T} h(\tau,\tilde{s})\right)T^{\alpha}\quad
\forall\tilde{s}\in(0,{\bar{s}}].
\end{aligned}
\end{equation}
Таким образом, в силу строгой монотонности функции \\$R_{\tilde{s}}(t):=\left(
E(t,\tilde{s})+\sup_{0<\tau<t} h(\tau,\tilde{s})\right)t^{\alpha}$ следует, что
$t_1(\tilde{s})<T$ $\forall\tilde{s}\in(0,{\bar{s}}]$. Отметим также, что в силу
\eqref{add1} из определения \eqref{11.5} вытекает:
\begin{equation}\label{11.6}
t'=t'(\tilde{s})<T\quad\forall\tilde{s}\in(0,\bar{s}].
\end{equation}
%\begin{equation}\label{11.7}
%t'=T\text{ если }\sup_{t\rightarrow T}\left( E(t,\tilde{s})+
%h(t,\tilde{s})\right)=\infty.
%\end{equation}
Итак, можем заключить, что функция $\Gamma_{\tilde{s}}(\cdot)$ определяет строго
монотонную возрастающую последовательность $\{t_j\}$ следующим соотношением:
\begin{equation}\label{11.9}
t_j:=\Gamma_{\tilde{s}}(t_{j-1}),\quad j=1,2,...,j_0=j_0(\tilde{s})<\infty:
t_{j_0}=\Gamma_{\tilde{s}}(t_{j_0-1})>t',\ \ \ t_{j_0-1}\leqslant t'.
\end{equation}
По последовательности $\{t_j\}$ из \eqref{11.9} определим последовательность интервалов
$\Delta_j=\Delta_j(\tilde{s})=t_j-t_{j-1}$, $j=1,2,...,j_0$, для которых в силу
определения \eqref{11.3} имеет место соотношение:
\begin{equation}\label{11.10}
\Delta_j^{-\alpha}=\frac{\xi^{\alpha}}{\omega_0T^{\alpha-\alpha_1}} \big(
E(t_j,\tilde{s})- E(t_{j-1},\tilde{s})+ \sup_{t_{j-1}<\tau<t_j} h(\tau,\tilde{s})\big).
\end{equation}
Покажем теперь, что последовательность $\{\Delta_j\}$ является квалифицированно
убывающей. Из условия теоремы \ref{Th1} следует, что для энергетических функций $E(t)$,
$h(t)$ рассматриваемого решения $u$ имеет место оценка:
\begin{equation}\label{add5}
E(t)+\sup_{0<\tau<t}h(\tau)=E_u(t)+\sup_{0<\tau<t}h_u(\tau)\leqslant\omega_0(T-t)^{-\alpha}\quad\forall\,t<T.
\end{equation}
В силу \eqref{add5} из определения \eqref{11.3} функции $\Gamma_{\tilde{s}}(t)$ вытекает
следующее неравенство:
\begin{equation*}
\begin{aligned}
&\Delta_j^{-\alpha}=\frac{\xi^{\alpha}}{\omega_0T^{\alpha-\alpha_1}} \big(
E(t_j,\tilde{s})- E(t_{j-1},\tilde{s})+ \sup_{t_{j-1}<\tau<t_j}
h(\tau,\tilde{s})\big)\leqslant
\\&\leqslant \frac{\xi^{\alpha}}
{\omega_0T^{\alpha-\alpha_1}} \left(E(t_j)+\sup_{0<\tau<t_j}h(\tau)\right)\leqslant
\xi^{\alpha}T^{-(\alpha-\alpha_1)}(T-t_j)^{-\alpha} \quad\forall j\leqslant j_0.
\end{aligned}
\end{equation*}
Отсюда в силу условия $T\geqslant1$ из формулировки задачи получаем:
\begin{equation}\label{4.6*}
\Delta_j\geqslant\xi^{-1}T^{1-\frac{\alpha_1}\alpha}(T-t_j)\geqslant
\xi^{-1}(T-t_j)\geqslant\xi^{-1}\Delta_{j+1}\quad\Rightarrow\quad
\Delta_{j+1}\leqslant\xi\Delta_j\quad\forall j\leqslant j_0,
\end{equation}
где $\Delta_{j_0+1}:=T-t_{j_0}$. Соотношение \eqref{4.6*} мы и называем свойством
квалифицированной монотонности последовательности $\{\Delta_j\}$. По сформированной
последовательности $\{t_j\}$ определим послойные энергетические функции $ E_j(s)$ и $
h_j(s)$ следующим образом:
\begin{equation}\label{4.7}
\begin{aligned}
 E_j(s):= E(t_j,{s})- E(t_{j-1},{s}),\quad
 h_j(s):=\sup_{t_{j-1}<\tau<t_j} h(\tau,{s})\quad\forall\,j\leqslant
j_0,\,\forall\,s\in(0,s_\Omega).
\end{aligned}
\end{equation}
В силу леммы \ref{Lem6.2.3} из \S3 энергетические функции из \eqref{4.7} удовлетворяют
системе \eqref{sys1}, \eqref{sys2}. Анализируя эту систему, установим оценки
энергетических функций $E(t,s)$, $h(t,s)$. Для этого сначала введем весовые
энергетические функции:
\begin{equation}\label{weight2}
A_j(s):=\Delta_j^{\alpha_1}E_j(s),\quad H_j(s):=\Delta_j^{\alpha_1}h_j(s),\quad
j=1,2,...,j_0.
\end{equation}
Для них стартовая система \eqref{sys1}, \eqref{sys2} принимает вид:
\begin{equation}\label{4.9}
\begin{aligned}
&  A_j(s)+H_j(s)\leqslant \overline{C}_1H_{j-1}(s)+C_2\Delta_j^{\nu_1-\alpha_1\mu_1}
    \left(-A'_j(s)\right)^{1+\mu_1}+
\\&C_3\Delta_j^{\nu_2-\alpha_1\mu_2}
    \left(-A'_j(s)\right)^{1+\mu_2}\quad\forall\,s\in(\tilde{s},\bar{s}),
    \\&H_j(s)\leqslant(1+\gamma)\xi^{\alpha_1}H_{j-1}(s)+ C_4
\gamma^{-(\nu_1+\mu_1)}
    \Delta_j^{\nu_1-\alpha_1\mu_1}
    \left(-A'_j(s)\right)^{1+\mu_1}+
\\&C_5 \gamma^{-\frac1q}
    \Delta_j^{\nu_2-\alpha_1\mu_2}
    \left(-A'_j(s)\right)^{1+\mu_2},
\end{aligned}
\end{equation}
где $\overline{C}_1=C_1\xi^{\alpha_1}$, $\overline{H}_0(s)=\Delta_0^{\alpha_1} h_0(s)$,
$\Delta_0=\xi^{-1}\Delta_1$. Наложим теперь первое условие на выбор постоянной $\xi$, а
именно положим $\xi<(1+\gamma)^{-\alpha_1^{-1}}$. При этом в силу \eqref{4.6*}:
\begin{equation}\label{4.10}
\lambda_j:=(1+\gamma)\left(\frac{\Delta_j}{\Delta_{j-1}}\right)^{\alpha_1}\leqslant
\lambda:=(1+\gamma)\xi^{\alpha_1}<1,\quad j=1,2, ...,j_0.
\end{equation}
Легко проверить следующее соотношение:
\begin{equation*}%\label{4.11}
\lambda_j\lambda_{j-1}...\lambda_{i+1}\Delta_i^{\nu_k-\alpha_1\mu_k}=
(1+\gamma)^{j-i}\Delta_j^{\nu_k-\alpha_1\mu_k}
\left(\frac{\Delta_j}{\Delta_{i}}\right)^{(1+\mu_k)\alpha_1-\nu_k}\quad\forall\,j>i,\,k=1,2.
\end{equation*}
Учитывая это, проитерируем неравенства \eqref{4.9}. В результате получаем:
\begin{equation}\label{4.12}
\begin{aligned}
&{A}_j(s)+{H}_j(s)\leqslant\overline{C}_1(1+\gamma)^{j-1}
\left(\frac{\Delta_j}{\Delta_0}\right)^{\alpha_1}{H}_0(s)+
\\&+\Delta_j^{\nu_1-\alpha_1\mu_1}C_6\gamma^{-(\nu_1+\mu_1)}\left[\sum_{i=1}^j(1+\gamma)^{j-i}
\left(\frac{\Delta_j}{\Delta_{i}}\right)^{(1+\mu_1)\alpha_1-\nu_1}\left(-{A}'_j(s)\right)^{1+\mu_1}\right]+
\\&+\Delta_j^{\nu_2-\alpha_1\mu_2}C_7\gamma^{-\frac1q}\left[\sum_{i=1}^j(1+\gamma)^{j-i}
\left(\frac{\Delta_j}{\Delta_{i}}\right)^{(1+\mu_2)\alpha_1-\nu_2}\left(-{A}'_j(s)\right)^{1+\mu_2}\right]
\\&\forall\,j\leqslant j_0,\,\forall\,s\in(\tilde{s},{\bar{s}}),
\end{aligned}
\end{equation}
где $C_6=\max\left\{C_2\gamma^{\left(\nu_1+\mu_1\right)},
\overline{C}_1C_4\lambda^{-1}\right\}$, $C_7=\max\left\{C_3\gamma^{\frac1q},
\overline{C}_1C_5\lambda^{-1}\right\}$. Преобразуем теперь систему \eqref{4.12} к такому
виду, чтобы к ней была применима лемма \ref{Lem9.2.6}. Для этого введем новые
энергетические функции
\begin{equation}\label{4.13}
\begin{aligned}
&{U}_j^{(1)}(s):=\sum_{i=1}^j(1+\gamma)^\frac{j-i}{1+\mu_1}
\left(\frac{\Delta_j}{\Delta_{i}}\right)^{\alpha_1-\frac{\nu_1}{1+\mu_1}}({A}_i(s)+{H}_i(s)),
\\&{U}_j^{(2)}(s):=\sum_{i=1}^j(1+\gamma)^\frac{j-i}{1+\mu_2}
\left(\frac{\Delta_j}{\Delta_{i}}\right)^{\alpha_1-\frac{\nu_2}{1+\mu_2}}({A}_i(s)+{H}_i(s)),
\quad j=1,2, ..., j_0.
\end{aligned}
\end{equation}
Очевидны соотношения:
\begin{equation}\label{4.14}
\begin{aligned}
&{U}_j^{(1)}(s)-A_j(s)-H_j(s)=\theta_{1,j}\,\overline{U}_{j-1}^{(1)}(s),\quad
{U}_0^{(1)}(s)=0,
\\&{U}_j^{(2)}(s)-A_j(s)-H_j(s)=\theta_{2,j}\,\overline{U}_{j-1}^{(2)}(s),\quad
{U}_0^{(2)}(s)=0,\quad j=1,2, ..., j_0,
\end{aligned}
\end{equation}
где
\begin{equation}\label{4.15}
\begin{aligned}
&\theta_{1,j}:=(1+\gamma)^\frac{1}{1+\mu_1}
\left(\frac{\Delta_j}{\Delta_{j-1}}\right)^{\alpha_1-\frac{\nu_1}{1+\mu_1}}\leqslant\theta_1:=
(1+\gamma)^\frac1{1+\mu_1}\xi^{\alpha_1-\frac{\nu_1}{1+\mu_1}}<1,
\\&\theta_{1,j}:=(1+\gamma)^\frac{1}{1+\mu_2}
\left(\frac{\Delta_j}{\Delta_{j-1}}\right)^{\alpha_1-\frac{\nu_2}{1+\mu_2}}\leqslant\theta_2:=
(1+\gamma)^\frac1{1+\mu_2}\xi^{\alpha_1-\frac{\nu_2}{1+\mu_2}}<1,
\end{aligned}
\end{equation}
Условиями \eqref{4.15} накладываются более жесткие окончательные требования на выбор
константы $\xi$. Очевидно, что $H_j(s)$, $j=1,2, ..., j_0$, являются абсолютно
непрерывными монотонно невозрастающими функциями. Поэтому из неравенств \eqref{4.12} в
силу соотношений \eqref{4.14} вытекает справедливость для почти всех
$s\in(\tilde{s},{\bar{s}})$ соотношений:
\begin{equation}\label{2.13}
\begin{aligned}
&U_j^{(1)}(s)\leqslant \overline{C}_1\lambda^{j-1}H_0(s)+\theta_1U_{j-1}^{(1)}(s)+
C_6\gamma^{-(\nu_1+\mu_1)}
\Delta_j^{\nu_1-\alpha_1\mu_1}\left(-\frac{d}{ds}U_j^{(1)}(s)\right)^{1+\mu_1}+
\\&+C_7\gamma^{-\frac1q}
\Delta_j^{\nu_2-\alpha_1\mu_2}\left(-\frac{d}{ds}U_j^{(2)}(s)\right)^{1+\mu_2}, \quad
j=1,2, ..., j_0,
\end{aligned}
\end{equation}
\begin{equation}\label{2.14}
\begin{aligned}
&U_j^{(2)}(s)\leqslant \overline{C}_1\lambda^{j-1}H_0(s)+\theta_2U_{j-1}^{(2)}(s)+
C_6\gamma^{-(\nu_1+\mu_1)}
\Delta_j^{\nu_1-\alpha_1\mu_1}\left(-\frac{d}{ds}U_j^{(1)}(s)\right)^{1+\mu_1}+
\\&+C_7\gamma^{-\frac1q}
\Delta_j^{\nu_2-\alpha_1\mu_2}\left(-\frac{d}{ds}U_j^{(2)}(s)\right)^{1+\mu_2}, \quad
j=1,2, ..., j_0.
\end{aligned}
\end{equation}
Оценим сверху значение ${U}_j^{(1)}(\tilde{s})$. В силу \eqref{weight2} и определения
\eqref{11.10} интервалов $\{\Delta_j\}=\{\Delta_j(\tilde{s})\}$ $\forall j\leqslant j_0$
имеем
\begin{equation}\label{4.16}
\begin{aligned}
&{U}_j^{(1)}(\tilde{s})=\sum_{i=1}^j(1+\gamma)^\frac{j-i}{1+\mu_1}
\left(\frac{\Delta_j}{\Delta_{i}}\right)^{\alpha_1-\frac{\nu_1}{1+\mu_1}}
\Delta_{i}^{\alpha_1}( E_i(\tilde{s})+ h_i(\tilde{s}))=
\\&=\omega_0\xi^{-\alpha}T^{\alpha-\alpha_1}\sum_{i=1}^j(1+\gamma)^\frac{j-i}{1+\mu_1}
\left(\frac{\Delta_j}{\Delta_{i}}\right)^{\alpha_1-\frac{\nu_1}{1+\mu_1}}
\Delta_{i}^{-(\alpha-\alpha_1)}\quad\forall\,j\leqslant j_0.
\end{aligned}
\end{equation}
Отсюда легко следует
\begin{equation}\label{4.16*}
\begin{aligned}
&{U}_j^{(1)}(\tilde{s})=\omega_0\xi^{-\alpha}T^{\alpha-\alpha_1}\Delta_{j}^{-(\alpha-\alpha_1)}
\sum_{i=1}^j(1+\gamma)^\frac{j-i}{1+\mu_1}
\left(\frac{\Delta_j}{\Delta_{i}}\right)^{\alpha-\frac{\nu_1}{1+\mu_1}}\leqslant
\\&\leqslant \omega_0\xi^{-\alpha}T^{\alpha-\alpha_1}\Delta_{j}^{-(\alpha-\alpha_1)}
\sum_{i=1}^j\left(\theta_1\xi^{\alpha-\alpha_1}\right)^{j-1}\leqslant
\\&\leqslant G_1\omega_0(1-\theta_1\xi^{\alpha-\alpha_1})^{-1}\Delta_{j}^{-(\alpha-\alpha_1)}
\quad\forall\,j\leqslant j_0,
\end{aligned}
\end{equation}
где $G_1=\xi^{-\alpha}T^{\alpha-\alpha_1}$. Аналогично получаем оценку для
$\overline{U}_j^{(2)}(\tilde{s})$:
\begin{equation}\label{4.17}
\overline{U}_j^{(2)}(\tilde{s})\leqslant
G_1\omega_0(1-\theta_2\xi^{\alpha-\alpha_1})^{-1}\Delta_{j}^{-(\alpha-\alpha_1)}
\quad\forall\,j\leqslant j_0.
\end{equation}
Складывая неравенства \eqref{2.13} и \eqref{2.14} и учитывая монотонное невозрастание
функций ${U}_j^{(1)}(s)$ и ${U}_j^{(2)}(s)$, получим следующее дифференциальное
неравенство относительно функций ${{U}}_j(s):={U}_j^{(1)}(s)+{U}_j^{(2)}(s)$:
\begin{equation}\label{2.16}
\begin{aligned}
&U_j(s)\leqslant
2\max\left\{\frac{C_6}{\gamma^{\nu_1+\mu_1}}\Delta_j^{\mu_1(\alpha+\beta-\alpha_1)}\left(-U'_j(s)\right)^{1+\mu_1},
\frac{C_7}{\gamma^{\frac1q}}\,\Delta_j^{\mu_2(\alpha+\beta-\alpha_1)}\left(-U'_j(s)\right)^{1+\mu_2}\right\}+
\\&+2\overline{C}_1\lambda^{j-1}H_0(s)+\overline{\theta}U_{j-1}(s)\quad\text{для почти всех
}s\in(\tilde{s},{\bar{s}}),\quad\overline{\theta}=\max(\theta_1,\theta_2).
\end{aligned}
\end{equation}
Соответственно, неравенства \eqref{4.17} и \eqref{4.16*} порождают "начальное" условие
для функций ${{U}}_j(s)$:
\begin{equation}\label{4.18}
{U}_j(\tilde{s})\leqslant G_2\omega_0\Delta_{j}^{-(\alpha-\alpha_1)}
\quad\forall\,j\leqslant j_0,
\end{equation}
где
$G_2=G_1\left((1-\theta_1\xi^{\alpha-\alpha_1})^{-1}+(1-\theta_2\xi^{\alpha-\alpha_1})^{-1}\right)$.
Несложно проверить, что систему \eqref{2.16}, \eqref{4.18} можно записать в следующем
"однородном" виде:
\begin{equation}\label{4.19}
\begin{aligned}
&\widetilde{U}_j(s):={U}_j(s)-\bar{b}\leqslant\overline{\theta}\widetilde{U}_{j-1}(s)+
\max\left\{k_j^{(1)}(-\widetilde{U}'_j(s))^{1+\mu_1},k_j^{(2)}(-\widetilde{U}'_j(s))^{1+\mu_2}\right\},
\\&\widetilde{U}_j(\tilde{s})\leqslant
K_j:=G_2\omega_0\Delta_{j}^{-(\alpha-\alpha_1)},\quad\bar{b}=2\overline{C}_1(1-\overline{\theta})^{-1}\overline{H}_0(0),\,
j=1,2, ... , j_0(\tilde{s}),
\end{aligned}
\end{equation}
где $k_j^{(1)}=\overline{C}_6\Delta_j^{\mu_1(\alpha+\beta-\alpha_1)}$,
$k_j^{(2)}=\overline{C}_7\Delta_j^{\mu_2(\alpha+\beta-\alpha_1)}$,
$\overline{C}_6=2C_6\gamma^{-(\nu_1+\mu_1)}$, $\overline{C}_7=2C_7\gamma^{-\frac1q}$. К
системе \eqref{4.19} уже применима лемма \ref{Lem9.2.6}. В силу этой леммы получим
следующую равномерную оценку:
\begin{equation}\label{4.20}
\begin{aligned}
&\widetilde{U}_j(s)\leqslant\omega_0^{\gamma_1}
\max\left\{G_3\psi(s-\tilde{s}),G_4\right\} \quad\forall\,s\in(\tilde{s},{\bar{s}}), \ \
\gamma_1:=\frac{\alpha+\beta-\alpha_1}{\beta},\ \ j=1,2, ... , j_0,
\end{aligned}
\end{equation}
где $\psi(s):=s^{{-\frac{(1+\mu_1)(\alpha-\alpha_1)}{\mu_1\beta}}}$,
\\$G_3=\left(\frac{\beta}{\alpha+\beta-\alpha_1}\right)^\frac{1+\mu_1}{\mu_1}
\left(\frac{\alpha-\alpha_1}{\alpha+\beta-\alpha_1}\right)^\frac{(1+\mu_1)(\alpha-\alpha_1)}{\mu_1\beta}
\bar{\mu}_1\bar{\mu}_2^\frac{(\alpha+\beta-\alpha_1)}{\beta}
\left(\frac{\overline{C}_6}{1-\overline{\theta}}\right)^\frac{\alpha-\alpha_1}{\mu_1\beta}
G_2^\frac{\alpha+\beta-\alpha_1}{\beta}$,
\\$G_4=\bar{\mu}_1\bar{\mu}_2
(1-\overline{\theta})^\frac{\alpha-\alpha_1}{\beta}
\overline{C}_6^{-\frac{(1+\mu_2)(\alpha-\alpha_1)}{(\mu_2-\mu_1)\beta}}
\overline{C}_7^{\frac{(1+\mu_1)(\alpha-\alpha_1)}{(\mu_2-\mu_1)\beta}}
G_2^\frac{\alpha+\beta-\alpha_1}{\beta}$,
$\bar{\mu}_1=\left(\frac{\mu_1}{1+\mu_1}\right)^\frac{1+\mu_1}{\mu_1} $,
\\$\bar{\mu}_2=\left(\frac{\mu_2}{1+\mu_2}\right)^\frac{(1+\mu_1)(1+\mu_2)}{\mu_1\mu_2}$.
Соответственно для $U_j(s)$ имеем оценку:
\begin{equation}\label{4.21}
U_j(s)\leqslant\omega_0^{\gamma_1} \max\left\{G_3\psi(s-\tilde{s}),G_4\right\}+\bar{b}
\quad\forall\,s\in(\tilde{s},{\bar{s}}),
\end{equation}
Далее определим значение $s_1>\tilde{s}$ следующим образом:
\begin{equation}\label{4.22}
\psi(s_1-\tilde{s})=G_3^{-1}\max\left\{G_4,\bar{b}\,\omega_0^{-\gamma_1}\right\}.
\end{equation}
Тогда из \eqref{4.21} вытекает справедливость неравенства:
\begin{equation}\label{4.23}
{{U}}_j(s)\leqslant2G_3\,\omega_0^{\gamma_1}\psi(s-\tilde{s})
\quad\forall\,s:\tilde{s}<s<s_2:=\min\{s_1,{\bar{s}}\},\quad\forall\,j\leqslant
j_0(\tilde{s}).
\end{equation}
Вспоминая определения \eqref{4.13} и \eqref{weight2}, выводим с помощью \eqref{4.23}
следующую оценку:
\begin{equation}\label{4.24}
 E_j(s)+ h_j(s)\leqslant2^{-1}\Delta_j^{-\alpha_1}{{U}}_j(s)\leqslant
G_3\omega_0^{\gamma_1}\psi(s-\tilde{s})\Delta_j^{-\alpha_1}
\quad\forall\,s\in(\tilde{s},s_2),\,j\leqslant j_0.
\end{equation}
Теперь оценим энергетические функции $E(t,s)$ и $h(t,s)$. Для этого зафиксируем
произвольное значение $i\leqslant j_0$ и просуммируем неравенства \eqref{4.24} по $j$ от
$1$ до $i$. В силу \eqref{11.10} и \eqref{4.6*} получаем
\begin{equation}\label{4.25}
\begin{aligned}
& E(t_i,s)+\sup_{0<\tau<t_i}h(\tau,s)\leqslant
G_3\omega_0^{\gamma_1}\psi(s-\tilde{s})\sum_{j=1}^i\Delta_j^{-\alpha_1}\leqslant
G_3\omega_0^{\gamma_1}\psi(s-\tilde{s})\Delta_i^{-\alpha_1}\sum_{j=1}^i\left(\xi^{\alpha_1}\right)^{j-1}\leqslant
\\&\leqslant G_3\omega_0^{\gamma_1}\psi(s-\tilde{s})\Delta_i^{-\alpha_1}(1-\xi^{\alpha_1})^{-1}=
G_5\omega_0^{\gamma_2}\psi(s-\tilde{s}) \left( E_i(\tilde{s})+
h_i(\tilde{s})\right)^\frac{\alpha_1}{\alpha} \quad\forall\,s\in(\tilde{s},s_2),
\\&\text{где }\gamma_2=\frac{(\alpha+\beta)(\alpha-\alpha_1)}{\alpha\beta},
\quad
G_5=\xi^{\alpha_1}(1-\xi^{\alpha_1})^{-1}T^{-\frac{(\alpha-\alpha_1)\alpha_1}{\alpha}}
G_3.
\end{aligned}
\end{equation}
Следующий шаг доказательства --- получение оценки типа \eqref{4.25} для произвольной
точки $t<T$. Для этого отметим, что функция $\Gamma_{\tilde{s}}(\cdot)$, определенная в
\eqref{11.3}, непрерывно, монотонно и взаимнооднозначно отображает любой отрезок
$[t_{j-1},t_j]$ на $[t_{j},t_{j+1}]$ $\forall\,j\leqslant j_0-1$. Зафиксируем
произвольную точку $\bar{t}\in[t_1,T)$. Пусть для определенности
$\bar{t}:=\bar{t}_k\in(t_k,t_{k+1}]$ при некотором $k\leqslant j_0$. Тогда единственным
образом восстановится последовательность $\{\bar{t}_i\}$, $i\leqslant k-1$ такая, что:
\begin{equation*}
\bar{t}_{i+1}=\Gamma_{\tilde{s}}(\bar{t}_{i})\quad\forall\,i\leqslant k-1,
\quad\bar{t}_i\in(t_i,t_{i+1}],\quad\bar{t}_0\in(0,t_{1}].
\end{equation*}
Этой последовательностью определяются новые смещения $\{\overline{\Delta}_i\}$:
\begin{equation}\label{4.26}
\overline{\Delta}_i^{-\alpha}:=(\bar{t}_{i}-\bar{t}_{i-1})^{-\alpha}=\frac{\xi^\alpha}{\omega_0T^{\alpha-\alpha_1}}
\Big( E(\bar{t}_{i},\tilde{s})- E(\bar{t}_{i-1},\tilde{s})+
\sup_{\bar{t}_{i-1}<t<\bar{t}_{i}} h(t,\tilde{s})\Big).
\end{equation}
Аналогично \eqref{4.6*} проверяется квалифицированная монотонность последовательности
$\{\overline{\Delta}_i\}:\overline{\Delta}_{i+1}<\xi\overline{\Delta}_i$
$\forall\,i\leqslant k-1$. По смещениям $\{\overline{\Delta}_i\}$ определяем
соответствующие энергетические функции $\overline{ E}_i(s)$ и $\overline{ h}_i(s)$.
Далее повторяем рассуждения \eqref{weight2}--\eqref{4.25}. Отличие будет состоять только
в том, что в качестве начальной функции выступает функция $\overline{ h}_0(s):=
h(\bar{t}_0,s)$. Оценим ${\overline{H}}_0(s)$, учитывая, что для энергетических функций
выполняется условие \eqref{setU} с функцией $F$, заданной соотношениями \eqref{LS-reg},
\eqref{regLS}:
\begin{equation*}
{\overline{H}}_0(s)=(\xi^{-1}\Delta_1)^{\alpha_1}\overline{ h}_0(s)\leqslant
\xi^{-\alpha_1}\omega_0\left(\frac{\bar{t}_1-\bar{t}_0}{T-\bar{t}_0}\right)^{\alpha_1}
\leqslant\xi^{-\alpha_1}\omega_0.
\end{equation*}
В результате получаем следующую оценку типа \eqref{4.25}:
\begin{equation}\label{4.27}
\begin{aligned}
& E(\bar{t},s)+ \sup_{0<\tau<\bar{t}}h(\tau,s)\leqslant
G_5\omega_0^{\gamma_2}\psi(s-\tilde{s}) \left( \overline{E}_k(\tilde{s})+
\overline{h}_k(\tilde{s})\right)^\frac{\alpha_1}{\alpha}
\quad\forall\,s\in(\tilde{s},s_3),
\end{aligned}
\end{equation}
где $\gamma_2$ и $G_5$ из \eqref{4.25}, $s_3:=\min\{\bar{s}_1,{\bar{s}}\}$, а значение
$\bar{s}_1$ определяется следующим соотношением:
\begin{equation*}
\psi(\bar{s}_1-\tilde{s})=G_3^{-1}\max\left\{G_4,\bar{\bar{b}}\,\omega_0^{-\gamma_1}\right\},
\quad\bar{\bar{b}}:=2\overline{C}_1(1-\overline{\theta})^{-1}{\overline{H}}_0(0).
\end{equation*}
Теперь в силу того, что $\bar{t}=\bar{t}_k$ является произвольной точкой из интервала
$(0,T]$, получаем из \eqref{4.27} оценку:
\begin{equation}\label{4.28}
\begin{aligned}
&E(t,s)+ \sup_{0<\tau<t}h(\tau,s)\leqslant
G_5\omega_0^{\frac{(\alpha+\beta)(\alpha-\alpha_1)}{\alpha\beta}}
(s-\tilde{s})^{-\frac{(1+\mu_1)(\alpha-\alpha_1)}{\mu_1\beta}} \left( E(t,\tilde{s})+
\sup_{0<\tau<t}h(\tau,\tilde{s})\right)^\frac{\alpha_1}{\alpha}
\\&\forall\,t\leqslant
T,\,\forall\,s,\,\tilde{s}:0<\tilde{s}<s<s_3,
\end{aligned}
\end{equation}
Таким образом установлено функциональное неравенство структуры \eqref{StampCond}
относительно параметрического семейства функций
$U_t(s):=E(t,s)+\sup_{0<\tau<t}h(\tau,s)$. В силу леммы \ref{StampLem} из \eqref{4.28}
вытекает следующая оценка:
\begin{equation}\label{4.29}
\begin{aligned}
& E(t,s)+ \sup_{0<\tau<t}h(\tau,s)\leqslant
2^{\frac{(1+\mu_1)\alpha^3}{\mu_1\alpha_1\beta(\alpha-\alpha_1)}}
G_5^{\frac{\alpha}{\alpha-\alpha_1}} \omega_0^{\frac{\alpha+\beta}{\beta}}
s^{-\frac{(1+\mu_1)\alpha}{\mu_1\beta}} \quad\forall\,t\leqslant
T,\,\forall\,s\in(0,s_3).
\end{aligned}
\end{equation}
Несложно проверить, что оценка \eqref{4.29} совпадает с доказываемой оценкой
\eqref{main} при $G=2^{\frac{(1+\mu_1)\alpha^3}{\mu_1\alpha_1\beta(\alpha-\alpha_1)}}
G_5^{\frac{\alpha}{\alpha-\alpha_1}}$, $\hat{s}=s_3$.

\vskip 5mm {\bf Доказательство следствия \ref{CorProfile}}

 Пусть $y$ --- произвольная
точка из $\Omega$ и
\begin{equation}\label{disty}
d(y)=dist(y,\partial\Omega):=s>0.
\end{equation}
При $p>1$, $q<1$ детальный анализ оценок максимума модуля обобщенных решений дважды
вырождающихся параболических уравнений из \cite{Ivan.max1}, \cite{Diben1}, \cite{SkrBur}
приводят с учетом \eqref{disty} к неравенству:
\begin{equation}\label{LocMaxEst}
\begin{aligned}
&u(T,y)\leqslant\gamma s^{-\frac{n+p+1}{(1+\lambda)q}}
\left(\int_{T-\xi}^T\int_{B_{\frac s2}(y)}|u(t,x)|^{p+\lambda
q}dx\,dt\right)^\frac1{(1+\lambda)q}+\gamma\left(\frac{s^{p+1}}\xi\right)^\frac 1{p-q} \
\ \forall\lambda\in(0,1],\\&\forall\,\xi\in(0,T),
\end{aligned}
\end{equation}
где $\gamma=\gamma(\lambda)\rightarrow\infty$ при $\lambda\rightarrow0$,
$B_r(y):=\{x:|x-y|<r\}$. При дополнительном условии $1<p<1+q$ можем считать $\lambda$
свободным параметром:
\begin{equation}\label{cond6}
0<\lambda\leqslant\frac{q+1-p}q.
\end{equation}
При этом $p+\lambda q\leqslant q+1$ и в силу неравенства Гельдера:
\begin{equation*}
\int_{B_{\frac s2}(y)}|u(t,x)|^{p+\lambda q}dx \leqslant c_0s^\frac{n(q+1-(p+\lambda
q))}{q+1}\left(\int_{B_{\frac s2}(y)}|u(t,x)|^{q+1}dx\right)^\frac{p+\lambda q}{q+1}.
\end{equation*}
Используя очевидное включение $B_\frac s2(y)\subset\Omega\left(\frac s2\right)$, а также
оценку \eqref{main}, продолжим последнее неравенство и получим:
\begin{equation*}
\int_{B_{\frac s2}(y)}|u(t,x)|^{p+\lambda q}dx \leqslant c_0s^{\nu_1}G^\frac{p+\lambda
q}{q+1} \omega_0^\frac{p+\lambda q}{\beta(p-q)},\quad \nu_1=\frac{n(q+1-(p+\lambda
q))}{q+1}-\nu\,\frac{p+\lambda q}{q+1}.
\end{equation*}
В силу этого неравенства оценка \eqref{LocMaxEst} приводит к:
\begin{equation}\label{LocEst2}
\begin{aligned}
&|u(T,y)|\leqslant\gamma c_0^\frac1{(1+\lambda)q}G^\frac{p+\lambda q}{(q+1)(1+\lambda)q}
\omega_0^\frac{p+\lambda q}{(p-q)(1+\lambda)q\beta}\xi^\frac1{(1+\lambda)q}s^{-\nu_2}
+\gamma\xi^{-\frac1{p-q}}s^\frac{p+1}{p-q},
\\&\nu_2=\frac{(n+\nu)(p+\lambda
q)}{(q+1)(1+\lambda)q}+\frac{p+1}{(1+\lambda)q}.
\end{aligned}
\end{equation}
Оптимизируя эту оценку по $\xi\in(0,T)$, получаем:
\begin{equation}\label{PointEst2}
\begin{aligned}
&|u(T,y)|\leqslant A\omega_0^{\frac1{\beta(p-q)}}s^{-\nu_3}\quad
\forall\,s<2\hat{s},\quad \nu_3=\frac{n(p-q)+(p+1)(q+1)-\beta(p-q)(p+1)}{\beta(p-q)^2},
\\&A=A(\lambda)=\gamma c_0^\frac1{p+\lambda q}G^\frac{1}{q+1} (p-q)^{-\frac{p-q}{p+\lambda
q}}(p+\lambda q)\big((1+\lambda)q\big)^{-\frac{(1+\lambda)q}{p+\lambda q}}.
\end{aligned}
\end{equation}
Оценка \eqref{PointEst2} очевидно эквивалентна доказываемой оценке \eqref{main1}
финального профиля решения $u$.

\section{\bf Приложения: вспомогательные утверждения}

В силу $C^2$--гладкости границы $\partial\Omega$ области $Q=(0,T)\times\Omega$, в
которой рассматривается задача \eqref{equat}--\eqref{initial}, существует число
$s_\Omega>0$ такое, что функция
\begin{equation}\label{s_Omega}
d(x):= dist(x,\partial\Omega)\in
C^2(\Omega\setminus\Omega(s))\quad\forall\,s\in(0,s_\Omega),
\end{equation}
где $\Omega(s)$ из \eqref{Omega(s)}. Соответственно, при этом $\partial\Omega(s)$
является $C^2$--гладким многообразием для всех $s\in(0,s_\Omega)$.

\begin{lemma}\label{Lem6.2.3}
Пусть $u(t,x)$ --- произвольное энергетическое решение задачи \eqref{equat} ---
\eqref{initial} из определения \ref{Def1}. И пусть промежуток $[0,T)$ разбит некоторой
монотонно возрастающей последовательностью точек $\{t_j\}$ ($j=1,2, ...,
j_0\leqslant\infty$, $t_0=0$) на промежутки $[t_{j-1},t_j)$ длиной
$\Delta_j:=t_{j}-t_{j-1}>0$. Тогда для почти всех $s\in(0,s_\Omega)$ послойные
энергетические функции рассматриваемого решения $u$
\begin{equation}\label{energy_j}
\begin{aligned}
E_j(s):=\int_{t_{j-1}}^{t_j}\int_{\Omega(s)}|\nabla_xu(t,x)|^{p+1}dxdt, \quad
h_j(s):=\sup_{t_{j-1}\leqslant t<t_j}\int_{\Omega(s)}|u(t,x)|^{q+1}dx
\end{aligned}
\end{equation}
удовлетворяет следующей системе дифференциальных неравенств:
\begin{equation}\label{sys1}
    E_j(s)+h_j(s)\leqslant C_1h_{j-1}(s)+C_2\Delta_j^{\nu_1}
    \left(-E'_j(s)\right)^{1+\mu_1}+C_3\Delta_j^{\nu_2}
    \left(-E'_j(s)\right)^{1+\mu_2},\ \ j=1,2, ..., j_0,
\end{equation}
\begin{equation}\label{sys2}
    h_j(s)\leqslant(1+\gamma)h_{j-1}(s)+ C_4\gamma^{-(\nu_1+\mu_1)}
    \Delta_j^{\nu_1}
    \left(-E'_j(s)\right)^{1+\mu_1}+C_5\gamma^{-\frac1q}\Delta_j^{\nu_2}
    \left(-E'_j(s)\right)^{1+\mu_2},
\end{equation}
где $\gamma>0$ --- произвольная константа, $C_i=const>0$ $\forall\,i=\overline{1,5}$
зависят только от известных параметров задачи и не зависят от $\omega_0$ и $\gamma$,
\begin{equation*}
\begin{aligned}
&\nu_1=\frac{(1-\theta)(q+1)}{q(p+1)+\theta(p-q)},\quad
\mu_1=\frac{(1-\theta)(p-q)}{q(p+1)+\theta(p-q)},
\\&\theta=\frac{n(p-q)+q+1}{n(p-q)+(q+1)(p+1)}<1,\quad \nu_2=\frac{(q+1)}{q(p+1)},\quad
\mu_2=\frac{(p-q)}{q(p+1)}.
\end{aligned}
\end{equation*}
\end{lemma}
Доказательство аналогично доказательству леммы 6.2.3 из \cite{KSSh}.

\begin{lemma}\label{Lem9.2.6}{\textbf{(см. Леммы 9.2.1--9.2.6 из \cite{KSSh})}}
Пусть некоторое семейство неотрицательных абсолютно непрерывных монотонно невозрастающих
функций $\{M_j(s)\}$, $j\leqslant j_0\leqslant\infty$, удовлетворяет для почти всех
$s\in(0,s_0), s_0>0$, системе дифференциальных неравенств:
\begin{equation}\label{sys9.2.6}
\begin{aligned}
&M_j(s)\leqslant\lambda M_{j-1}(s)+(1-\lambda)
\max\left\{k_j^{(1)}(-M'_j(s))^{1+\gamma_1};k_j^{(2)}(-M'_j(s))^{1+\gamma_2}\right\}\quad\forall\,s>0,
\\&M_j(0)\leqslant K_j\quad\forall\,j\in\mathbb{N},\quad M_0(s):=0,
\end{aligned}
\end{equation}
где $\gamma_2>\gamma_1>0$, $\lambda=const\in(0,1)$,
$k_j^{(1)}=c_1\varepsilon_j^{\gamma_1}$, $k_j^{(2)}=c_2\varepsilon_j^{\gamma_2}$,
$K_j=c_3\varepsilon_j^{-(1-\delta)}$, $c_1,\,c_2,\,c_3>0$ --- произвольные постоянные,
$\delta\in(0,1)$, $\{\varepsilon_j\}$ --- произвольная монотонно убывающая
последовательность положительных чисел. Тогда для $M_j(s)$ справедлива следующая
равномерная оценка:
\begin{equation}\label{stat9.2.6}
M_j(s)\leqslant
\max\left\{B_1s^{-\frac{(1+\gamma_1)(1-\delta)}{\delta\gamma_1}},B_2\right\}
\quad\forall\,s\in(0,s_0),\,\forall\,j\leqslant j_0,
\end{equation}
где $B_1,\,B_2$ --- положительные константы, зависящие от
$c_1,\,c_2,\,c_3,\,\delta,\,\gamma_1,\,\gamma_2$, но не зависящие ни от $j_0$, ни от
последовательности $\{\varepsilon_j\}$.
\end{lemma}

\begin{lemma}\label{StampLem}{\textbf{(см. \cite{Stamp})}}
Пусть некоторая непрерывная неотрицательная невозрастающая функция
$f:(0,s_0)\rightarrow\mathbb{R}^1_+$ удовлетворяет соотношению:
\begin{equation}\label{StampCond}
f(s+\delta)\leqslant
a\delta^{-\rho}f(s)^\lambda\quad\forall\,s>0,\,\delta>0:s+\delta<s_0,\,0<\lambda=const<1,
\end{equation}
где $0<a=const<\infty$, $0<\rho=const<\infty$. Тогда для функции $f$ справедлива
следующая универсальная априорная оценка:
$f(s)\leqslant2^{\frac\rho{\lambda(1-\lambda)^2}}
a^{\frac1{1-\lambda}}s^{-\frac\rho{1-\lambda}}$ \ \ $\forall\,s\in(0,s_0)$.
\end{lemma}

\end{document}